\newfont{\Bbb}{msbm10 scaled\magstephalf}
\documentclass[reqno]{amsart}
\usepackage{amssymb}
\usepackage{amsfonts}
\usepackage{amsmath, amsthm, amscd, amsfonts}
\usepackage{cases}
\usepackage{cite}

\numberwithin{equation}{section}

\begin{document}
\title{Essential Norms of difference of generalized composition Operators from $\alpha$-Bloch spaces to $\beta$-Bloch spaces}

\author[N. Xu and Z.H. Zhou]  {Ning Xu and Ze-Hua Zhou$^*$}

\address{\newline Ning Xu\newline School of science, Jiangsu ocean University\newline Lianyungang 222005\newline P.R. China}
\email{gx899200@126.com}

\address{\newline  Ze-Hua Zhou\newline  School of Mathematics, Tianjin University, Tianjin 300354, P.R. China.\newline
} \email{zehuazhoumath@aliyun.com;
zhzhou@tju.edu.cn}



\date{}
\thanks{\noindent $^*$Corresponding author.\\
The work was supported in part by the National Natural Science Foundation of
China(Grant No.11771323) and the Scientific Research Foundation for Ph.D.of Jiangsu Ocean University(No.KQ17006).}

\begin{abstract}
In this paper, we study the boundedness and essential norms of the differences of two generalized composition operators acting from $\alpha$-Bloch space to $\beta$-Bloch space on the open unit disk. From essential norms, we get the compactness of the differences of two generalized composition operators. This study has a relationship to the topological structure of generalized composition operators acting from $\alpha$-Bloch space to $\beta$-Bloch space.

\end{abstract}

\maketitle

\section{Introduction}
Let $\mathbb{D}$ be an open disk in the complex plane $\mathbb{C}$ and $H(\mathbb{D})$ be the class of all functions analytic in $\mathbb{D}$.
We denote by $S(\mathbb{D})$ the set of all analytic self-maps on $\mathbb{D}$. For $0<\alpha<\infty$, recall that the Bloch type space $\mathcal{B}^\alpha$, or $\alpha$-Bloch space, consists of all $f\in H(\mathbb{D})$ such that
\begin{eqnarray*}
\|f\|_{\alpha}=\sup\limits_{z\in\mathbb{D}}(1-|z|^2)^\alpha|f'(z)|<\infty.
\end{eqnarray*}
It is well known that $\mathcal{B}^\alpha$ is a Banach space under the norm $\|f\|_{\mathcal{B}^\alpha}=|f(0)|+\|f\|_\alpha$. The little $\alpha$-Bloch $\mathcal{B}_0^\alpha$
is a subspace of $\mathcal{B}^\alpha$, consisting of all $f\in H(\mathbb{D})$ such that $\lim\limits_{|z|\rightarrow 1}(1-|z|^2)^\alpha|f'(z)|=0$.
When $\alpha=1$, $\mathcal{B}^\alpha$ is classical space $\mathcal{B}$.
\

For $z,w\in \mathbb{D}$, the pseudo-hyperbolic distance $\rho(z,w)$ between $z$ and $w$ is given by
\begin {eqnarray*}
\rho(z,w)=\Big|\frac{w-z}{1-\overline{w}z}\Big|.
\end{eqnarray*}
For $\lambda\in \mathbb{D}$, Let $\sigma_\lambda$ be the M\"{o}bius transformation of $\mathbb{D}$ defines by
\begin{eqnarray*}
\sigma_\lambda(z)=\frac{\lambda-z}{1-\overline{\lambda}z}.
\end{eqnarray*}
We remark that $\rho(\lambda,z)=|\sigma_\lambda(z)|\leq 1$.

An analytic self-map $\varphi\in S(\mathbb{D})$ induces the composition operator $C_\varphi$, defined by
$C_\varphi(f)=f(\varphi(z))$
for $f$ analytic on $S(\mathbb{D})$. The composition operator has been studied by many researchers on various spaces (see, for example,\cite{CB} and references therein). Motivated by the fact that composition operators and weighted composition operators naturally come from isometries of some functions spaces, in \cite{LS1}, for $g\in H(D)$ and $\varphi\in S(\mathbb{D})$, Li and Stevi\'{c} defined the generalized composition operator$ C_\varphi^g$ as follows:
\begin{eqnarray*}
C^g_\varphi f(z)=\int_0^zf'(\varphi(\xi))g(\xi)d\xi\ \ z\in \mathbb{D}.
\end{eqnarray*}
The boundedness and compactness of the generalized composition operator on Zygmund spaces and Bloch-type spaces were investigated in \cite{LS1}. Some related results concerning the generalized composition operator on various spaces can be found in, for example, \cite{L1,P,S1,ZL}.

The essential norm of a bounded linear operator $T:X\rightarrow Y$ is its distance to the set of compact operators $K$ mapping $X$ to $Y$, that is,
\begin{eqnarray*}
\|T\|_{e,X\rightarrow Y}=\inf\{\|T-K\|_{X\rightarrow Y}: K \ \rm{is  \ compact}\},
\end{eqnarray*}
where $X$, $Y$ are Banach spaces and $\|\cdot\|_{X\rightarrow Y}$ is the operator nom.

Recently, many researchers have been studying the difference of two operators. The study of the difference of two operators is motivated by the difference of two composition operators  which was started on $H^2$. The main purpose for this study is to understand the topological structure of $\mathcal{C}(H^2)$, the set of composition operators on $H^2$, see \cite{SS}. Then, MacCluer, Ohno and Zhao \cite{MOZ} considered the above problems on H$^\infty$. These works gave a relationship between a component problem and the behavior of the difference of two composition operators acting from $\mathcal{B}$ to H$^\infty$. After that, the study of the difference of two operators, such as boundedness, compactness and essential norms have been studied on several spaces of analytic functions by many authors, see \cite{HO1,FZ,YZ,SL,HL,KSK,CZ}.

As an analogue, the topological structure of the set $\mathcal{C}(\mathcal{B}^\alpha\rightarrow\mathcal{B}^\beta)$ of bounded generalized composition operators from $\mathcal{B}^\alpha$ to $\mathcal{B}^\beta$ can be considered. In this context, we deal with the differences of generalized composition operators from $\mathcal{B}^\alpha$ to $\mathcal{B}^\beta$. The main purpose of this paper is to express the boundedness and essential norms of  $C_\varphi^g-C_\psi^h$ from $\mathcal{B}^\alpha$ to $\mathcal{B}^\beta$ which generalize \cite{L2}. The authors expect that these results will play some roles in the succeeding investigation.

For two quantities $A$ and $B$ which may depend on $\varphi$ and $\psi$, we use the abbreviation $A\lesssim B$ whenever this is a positive constant $C$(independent of $\varphi$ and $\psi$) such that $A\leq CB$. We write $A\sim B$ if $A\lesssim B \lesssim A$.

\section{Prerequisites}
In this section, we give some notions and auxiliary lemmas needed in the proof of the following theorems.
We use the following notion:
\begin{eqnarray*}
\rho(z)=\rho(\varphi(z),\psi(z)),\ \ \
\tau(z)=\frac{(1-|\varphi(z)|^2)(1-|\psi(z)|^2)}{(1-\overline{\varphi(z)}\psi(z))^2}.
\end{eqnarray*}
We remark that
\begin{eqnarray*}
|\sigma_{\varphi(z)}(\psi(z))|=\rho(z)
\end{eqnarray*}
and $|\tau(z)|=1-\rho^2(z)$.

For convenience for the statements of our main results, we note
\begin{eqnarray*}
D_{\varphi,g}(z):=\frac{(1-|z|^\beta)g(z)}{(1-|\varphi(z)|^2)^\alpha}\ \ \ D_{\psi,h}(z):=\frac{(1-|z|^\beta)h(z)}{(1-|\psi(z)|^2)^\alpha}.
\end{eqnarray*}

{\bf Lemma 2.1.{\cite{SL}}} Let $0<\alpha<\infty$. For all $z,w\in \mathbb{D}$, the Bloch-type induced distance is given by
\begin{eqnarray*}
\flat_\alpha(z,w)=\sup\limits_{\|f\|_{\mathcal{B}^\alpha}\leq 1}|(1-|z|^2)^\alpha f'(z)-(1-|w|^2)^\alpha f'(w)|.
\end{eqnarray*}
Then
\begin{eqnarray*}
\flat_\alpha(z,w)\lesssim \rho(z,w).
\end{eqnarray*}

{\bf Remark 2.1.} Especially, define the function $f_a(z)=\int^z_0\frac{(1-|a|^2)^\alpha}{(1-\overline{a}u)^{2\alpha}}du$ see \cite{SL}. Then $\|f_a\|_{\mathcal{B}^\alpha}=1$. From Lemma 2.1, take $a=\varphi(z)$, we have that
\begin{eqnarray*}
|1-\tau^{\alpha}(z)|&=&\Big|1-\frac{(1-|\varphi(z)|^2)^{\alpha}(1-|\psi(z)|^2)^{\alpha}}{(1-\overline{\varphi(z)}\psi(z))^{2\alpha}}\Big|\\
&=&|(1-|\varphi(z)|^2)^{\alpha}f'_{\varphi(z)}(\varphi(z))-(1-|\psi(z)|^2)^{\alpha}f'_{\varphi(z)}(\psi(z))|\\
&\leq&\flat_\alpha(\varphi(z),\psi(z))\lesssim \rho(z).
\end{eqnarray*}

The following lemma can be proved in a standard way(see,e.g.,Proposition 3.11 in \cite{CB})

{\bf Lemma 2.2} Assume that $\varphi,\psi\in S(\mathbb{D}),$ $g, h\in S(\mathbb{D})$. Then $C^g_\varphi-C^h_\psi:\mathcal{B}^\alpha\rightarrow\mathcal{B}^\beta$ is compact if and only if $C^g_\varphi-C^h_\psi$ is bounded and for any bounded sequence $\{f_k\}_{k\in\mathbb{N}}$ in $\mathcal{B}^\alpha$ which converges to 0 uniformly on compact subsets of $\mathbb{D}$, $\|(C^g_\varphi-C^h_\psi)f_k\|_{\mathcal{B}^\beta}\rightarrow0$ as $k\rightarrow\infty$.

\section{Boundedness of $C^g_\varphi-C^h_\psi$: $\mathcal{B}^\alpha\rightarrow\mathcal{B}^\beta$}
In this section we give the characterization for the boundedness of the differences of generalized composition operators from $\mathcal{B}^\alpha$ to $\mathcal{B}^\beta$. For any $a\in \mathbb{D}$, we define the following test functions see \cite{SL}
\begin{eqnarray*}
g_a(z)=\int_0^z\frac{(1-|a|^2)^\alpha}{(1-\overline{a}u)^{2\alpha}}\frac{a-u}{1-\overline{a}u}du.
\end{eqnarray*}
It is easy to see that $\|g_a\|\leq\|f_a\|_{\mathcal{B}^\alpha}=1$. To prove the result in this section, we need the following Propositions.

{\bf Proposition 3.1.} Let $0<\alpha,\beta<\infty$. Let $\varphi, \psi\in S(\mathbb{D})$. Then the following inequalities hold:

(i) $\sup\limits_{z\in\mathbb{D}}|D_{\varphi,g}(z)|\rho(z)\lesssim
\sup\limits_{a\in\mathbb{D}}\|(C^g_\varphi-C^h_\psi)f_a\|_{\mathcal{B}^\beta}+\sup\limits_{a\in\mathbb{D}}\|(C^g_\varphi-C^h_\psi)g_a\|_{\mathcal{B}^\beta}.$

(ii)$\sup\limits_{z\in\mathbb{D}}|D_{\psi,h}(z)|\rho(z)\lesssim
\sup\limits_{a\in\mathbb{D}}\|(C^g_\varphi-C^h_\psi)f_a\|_{\mathcal{B}^\beta}+\sup\limits_{a\in\mathbb{D}}\|(C^g_\varphi-C^h_\psi)g_a\|_{\mathcal{B}^\beta}.$

(iii)$\sup\limits_{z\in\mathbb{D}}|D_{\varphi,g}(z)-D_{\psi,h}(z)|\lesssim
\sup\limits_{a\in\mathbb{D}}\|(C^g_\varphi-C^h_\psi)f_a\|_{\mathcal{B}^\beta}+\sup\limits_{a\in\mathbb{D}}\|(C^g_\varphi-C^h_\psi)g_a\|_{\mathcal{B}^\beta}.$

{\bf Proof.} For any $z\in \mathbb{D}$, take $a=\varphi(z)$, we have that
\begin{eqnarray*}
\|(C^g_\varphi-C^h_\psi)f_{\varphi(z)}\|_{\mathcal{B}^\beta}
&\geq&\Big|\frac{(1-|z|^2)^\beta}{(1-|\varphi(z)|^2)^\alpha}g(z)
-\frac{(1-|z|^2)^\beta(1-|\varphi(z)|^2)^\alpha}{(1-\overline{\varphi(z)}\psi(z))^{2\alpha}}h(z)\Big|\\
&\geq&|D_{\varphi,g}(z)|-|\tau^\alpha(z)||D_{\psi,h}(z)|,
\end{eqnarray*}
and
\begin{eqnarray*}
\|(C^g_\varphi-C^h_\psi)g_{\varphi(z)}\|_{\mathcal{B}^\beta}
&\geq&(1-|z|^2)^\beta\Big|g_{\varphi(z)}'(\varphi(z))g(z)-g_{\varphi(z)}'(\psi(z))h(z)\Big|\\
&=&|\tau^\alpha(z)||D_{\psi,h}(z)|\rho(z).
\end{eqnarray*}
Thus
\begin{eqnarray*}
|D_{\varphi,g}(z)|\rho(z)
&\leq&\|(C^g_\varphi-C^h_\psi)f_{\varphi(z)}\|_{\mathcal{B}^\beta}\rho(z)+\|(C^g_\varphi-C^h_\psi)g_{\varphi(z)}\|_{\mathcal{B}^\beta}\\
&\leq&\|(C^g_\varphi-C^h_\psi)f_{\varphi(z)}\|_{\mathcal{B}^\beta}+\|(C^g_\varphi-C^h_\psi)g_{\varphi(z)}\|_{\mathcal{B}^\beta}.
\end{eqnarray*}
Similarly
\begin{eqnarray*}
|D_{\psi,h}(z)|\rho(z)
&\leq&\|(C^g_\varphi-C^h_\psi)f_{\psi(z)}\|_{\mathcal{B}^\beta}+\|(C^g_\varphi-C^h_\psi)g_{\psi(z)}\|_{\mathcal{B}^\beta}.
\end{eqnarray*}
Hence, we have
\begin{eqnarray*}
|D_{\varphi,g}(z)|\rho(z)
&\leq&\|(C^g_\varphi-C^h_\psi)f_a\|_{\mathcal{B}^\beta}+\|(C^g_\varphi-C^h_\psi)g_a\|_{\mathcal{B}^\beta}.
\end{eqnarray*}
Similarly
\begin{eqnarray*}
|D_{\psi,h}|\rho(z)
&\leq&\|(C^g_\varphi-C^h_\psi)f_a\|_{\mathcal{B}^\beta}+\|(C^g_\varphi-C^h_\psi)g_a\|_{\mathcal{B}^\beta}.
\end{eqnarray*}
By Lemma 2.1, we obtain that
\begin{eqnarray*}
\|(C^g_\varphi-C^h_\psi)f_{\varphi(z)}\|_{\mathcal{B}^\beta}
&\geq&|D_{\varphi,g}(z)-D_{\psi,h}(z)|-|1-\tau^\alpha(z)||D_{\psi,h}(z)|\\
&\gtrsim&|D_{\varphi,g}(z)-D_{\psi,h}(z)|-|D_{\psi,h}(z)|\flat_\alpha(\varphi(z),\psi(z))\\
&\gtrsim&|D_{\varphi,g}(z)-D_{\psi,h}(z)|-|D_{\psi,h}(z)|\rho(z).
\end{eqnarray*}
Thus
\begin{eqnarray*}
|D_{\varphi,g}(z)-D_{\psi,h}(z)|&\leq&\|(C^g_\varphi-C^h_\psi)f_{\varphi(z)}\|_{\mathcal{B}^\beta}+|D_{\psi,h}(z)|\rho(z)\\
&\lesssim&\|(C^g_\varphi-C^h_\psi)f_{\varphi(z)}\|_{\mathcal{B}^\beta}+\|(C^g_\varphi-C^h_\psi)f_{\psi(z)}\|_{\mathcal{B}^\beta}\\
&+&\|(C^g_\varphi-C^h_\psi)g_{\psi(z)}\|_{\mathcal{B}^\beta}\\
&\lesssim&\|(C^g_\varphi-C^h_\psi)f_a\|_{\mathcal{B}^\beta}+\|(C^g_\varphi-C^h_\psi)g_a\|_{\mathcal{B}^\beta},
\end{eqnarray*}
which implies that (iii) holds.

{\bf Proposition 3.2.} Let $0<\alpha,\beta<\infty$. Let $\varphi, \psi\in S(\mathbb{D})$. Then the following inequalities hold:

(i) $\sup\limits_{a\in\mathbb{D}}\|(C^g_\varphi-C^h_\psi)f_a\|_{\mathcal{B}^\beta}
\leq\sup\limits_{n\in\mathbb{N}}\sup\limits_{z\in\mathbb{D}}(n+1)^\alpha(1-|z|^2)^\beta\Big|\varphi^n(z)g(z)-\psi^n(z)h(z)\Big|.$

(ii) $\sup\limits_{a\in\mathbb{D}}\|(C^g_\varphi-C^h_\psi)g_a\|_{\mathcal{B}^\beta}
\leq\sup\limits_{n\in\mathbb{N}}\sup\limits_{z\in\mathbb{D}}(n+1)^\alpha(1-|z|^2)^\beta\Big|\varphi^n(z)g(z)-\psi^n(z)h(z)\Big|.$

{\bf Proof.} For any $a\in\mathbb{D}$, recall that,
\begin{eqnarray*}
\frac{1}{(1-\overline{a}u)^{2\alpha}}=\sum_{k=0}^{\infty}\frac{\Gamma(k+2\alpha)}{k!\Gamma(2\alpha)}\overline{a}^ku^k,\ \ u\in\mathbb{D}.
\end{eqnarray*}
Then, we express $f_a$ into expansion as
\begin{eqnarray*}
f_a(z)=(1-|a|^2)^\alpha\int_0^z\sum_{k=0}^{\infty}\frac{\Gamma(k+2\alpha)}{k!\Gamma(2\alpha)}\overline{a}^ku^kdu.
\end{eqnarray*}
Noticing that $((C^g_\varphi-C^h_\psi)f_a)(0)=0$, hence
\begin{eqnarray*}
\|(C^g_\varphi-C^h_\psi)f_a\|_{\mathcal{B}^\beta}
&\leq&\sup\limits_{z\in\mathbb{D}}(1-|z|^2)^\beta(1-|a|^2)^\alpha\sum_{k=0}^{\infty}\frac{\Gamma(k+2\alpha)}{k!\Gamma(2\alpha)}|a|^k\\
&\times&\Big|\varphi^k(z)g(z)-\psi^k(z)h(z)\Big|\\
&\leq&(1-|a|^2)^\alpha\sum_{k=0}^{\infty}\frac{\Gamma(k+2\alpha)}{k!\Gamma(2\alpha)}(k+1)^{-\alpha}|a|^k\\
&\times&\sup\limits_{z\in\mathbb{D}}(k+1)^\alpha(1-|z|^2)^\beta\Big|\varphi^k(z)g(z)-\psi^k(z)h(z)\Big|,
\end{eqnarray*}
by Stirling's formula, $\frac{\Gamma(k+2\alpha)}{k!\Gamma(2\alpha)}(k+1)^{-\alpha}\approx k^{\alpha-1}, k\rightarrow\infty$. Note that (see \cite{Z}).
\begin{eqnarray*}
\frac{1}{(1-|a|^2)^\alpha}=\sum_{k=0}^{\infty}\frac{\Gamma(k+\alpha)}{k!\Gamma(\alpha)}|a|^k \ \ and \ \
\frac{\Gamma(k+\alpha)}{k!}\approx k^{\alpha-1},\ k\rightarrow \infty,
\end{eqnarray*}
we have
\begin{eqnarray*}
\sum_{k=0}^{\infty}\frac{\Gamma(k+2\alpha)}{k!\Gamma(2\alpha)}(k+1)^{-\alpha}|a|^k\approx \sum_{k=0}^{\infty}(k+1)^{\alpha-1}|a|^k\approx \frac{1}{(1-|a|^2)^\alpha},
\end{eqnarray*}
thus we have
\begin{eqnarray*}
\|(C^g_\varphi-C^h_\psi)f_a\|_{\mathcal{B}^\beta}
\lesssim\sup\limits_{n\in\mathbb{N}}\sup\limits_{z\in\mathbb{D}}(n+1)^\alpha(1-|z|^2)^\beta\Big|\varphi^n(z)g(z)-\psi^n(z)h(z)\Big|.
\end{eqnarray*}
Therefore (i) holds.

Note that
\begin{eqnarray*}
\frac{a-u}{1-\overline{a}u}=a-(1-|a|^2)\sum_{k=0}^{\infty}\overline{a}^ku^{k+1}.
\end{eqnarray*}
By an analogous calculation, we have that (see \cite{SL})
\begin{eqnarray*}
g_a(z)=af_a(z)-(1-|a|^2)^{\alpha+1}\int_0^z\sum_{k=1}^{\infty}\Big(\sum_{l=0}^{k-1}\frac{\Gamma(l+2\alpha)}{l!\Gamma(2\alpha)}\Big)\overline{a}^{k-1}u^kdu
\end{eqnarray*}
Therefore,
\begin{eqnarray*}
& &\|(C^g_\varphi-C^h_\psi)g_a\|_{\mathcal{B}^\beta}\leq\|(C^g_\varphi-C^h_\psi)f_a\|_{\mathcal{B}^\beta}\\
&+&\sup\limits_{z\in\mathbb{D}}(1-|a|^2)^{\alpha+1}(1-|z|^2)^\beta\sum_{k=1}^{\infty}\Big(\sum_{l=0}^{k-1}\frac{\Gamma(l+2\alpha)}{l!\Gamma(2\alpha)}\Big)
|a|^{k-1}\\
&\times&\Big|\varphi^k(z)g(z)-\psi^k(z)h(z)\Big|\\
&=&\|(C^g_\varphi-C^h_\psi)f_a\|_{\mathcal{B}^\beta}+(1-|a|^2)^{\alpha+1}\sum_{k=1}^{\infty}
\Big(\sum_{l=0}^{k-1}\frac{\Gamma(l+2\alpha)}{l!\Gamma(2\alpha)}\Big)|a|^{k-1}\\
&\times&\sup\limits_{z\in\mathbb{D}}(1-|z|^2)^\beta\Big|\varphi^k(z)g(z)-\psi^k(z)h(z)\Big|,
\end{eqnarray*}
by Stirling's formula,
\begin{eqnarray*}
\sum_{l=0}^{k-1}\frac{\Gamma(l+2\alpha)}{l!\Gamma(2\alpha)}\approx \sum_{l=0}^{k-1}l^{2\alpha}\approx k^{2\alpha},\ k\rightarrow\infty.
\end{eqnarray*}
Hence
\begin{eqnarray*}
& &(1-|a|^2)^{\alpha+1}\sum_{k=1}^{\infty}
\Big(\sum_{l=0}^{k-1}\frac{\Gamma(l+2\alpha)}{l!\Gamma(2\alpha)}\Big)|a|^{k-1}\\
&\times&\sup\limits_{z\in\mathbb{D}}(1-|z|^2)^\beta\Big|\varphi^k(z)g(z)-\psi^k(z)h(z)\Big|\\
&\leq&(1-|a|^2)^{\alpha+1}\sum_{k=1}^{\infty}
\Big(\sum_{l=0}^{k-1}\frac{\Gamma(l+2\alpha)}{l!\Gamma(2\alpha)}\Big)|a|^{k-1}k^{-\alpha}\\
&\times&\sup\limits_{n\in\mathbb{N}}\sup\limits_{z\in\mathbb{D}}n^\alpha(1-|z|^2)^\beta\Big|\varphi^n(z)g(z)-\psi^n(z)h(z)\Big|\\
&\lesssim&(1-|a|^2)^{\alpha+1}\sum_{k=1}^{\infty}k^\alpha|a|^{k-1}
\sup\limits_{n\in\mathbb{N}}\sup\limits_{z\in\mathbb{D}}n^\alpha(1-|z|^2)^\beta\Big|\varphi^n(z)g(z)-\psi^n(z)h(z)\Big|\\
&\lesssim&\sup\limits_{n\in\mathbb{N}}\sup\limits_{z\in\mathbb{D}}(n+1)^\alpha(1-|z|^2)^\beta\Big|\varphi^n(z)g(z)-\psi^n(z)h(z)\Big|.
\end{eqnarray*}
Thus
\begin{eqnarray*}
\|(C^g_\varphi-C^h_\psi)g_a\|_{\mathcal{B}^\beta}\lesssim
\sup\limits_{n\in\mathbb{N}}\sup\limits_{z\in\mathbb{D}}(n+1)^\alpha(1-|z|^2)^\beta\Big|\varphi^n(z)g(z)-\psi^n(z)h(z)\Big|.
\end{eqnarray*}
Hence (ii) holds.

{\bf Theorem 3.3.}
Let $0<\alpha,\beta<\infty$. Let $\varphi,\psi\in S(\mathbb{D})$.
Then the following statements are equivalent.

(i) $C^g_\varphi-C^h_\psi:\mathcal{B}^\alpha\rightarrow\mathcal{B}^\beta$ is bounded.

(iia)$\sup\limits_{z\in\mathbb{D}}|D_{\varphi,g}(z)-D_{\psi,h}(z)|+\sup\limits_{z\in\mathbb{D}}|D_{\psi,h}(z)|\rho(z)<\infty.$

(iib)$\sup\limits_{z\in\mathbb{D}}|D_{\varphi,g}(z)-D_{\psi,h}(z)|+\sup\limits_{z\in\mathbb{D}}|D_{\varphi,g}(z)|\rho(z)<\infty.$

(iii)$\|(C^g_\varphi-C^h_\psi)f_a\|_{\mathcal{B}^\beta}+\|(C^g_\varphi-C^h_\psi)g_a\|_{\mathcal{B}^\beta}<\infty.$

(iv)$\sup\limits_{n\in\mathbb{N}}\sup\limits_{z\in\mathbb{D}}(n+1)^\alpha(1-|z|^2)^\beta|\varphi^n(z)g(z)-\psi^n(z)h(z)|<\infty.$

{\bf Proof.} (i)$\Rightarrow$(iv) Suppose that $C^g_\varphi-C^h_\psi:\mathcal{B}^\alpha\rightarrow\mathcal{B}^\beta$ is bounded.
When $n\geq1$. Consider the function $z^n$. From \cite{Zh}, we see that $\|z^n\|_{\mathcal{B}^\alpha}\lesssim n^{1-\alpha}$. Let $f_n(z)=z^n/\|z^n\|_{\mathcal{B}^\alpha}$, then $\|f_n\|_{\mathcal{B}^\alpha}=1$. Hence
\begin{eqnarray*}
\infty&>&\|C^g_\varphi-C^h_\psi\|_{\mathcal{B}^\alpha\rightarrow\mathcal{B}^\beta}\geq\|(C^g_\varphi-C^h_\psi)f_n\|_{\mathcal{B}^\beta}\\
&\gtrsim&\sup\limits_{z\in\mathbb{D}}\frac{n}{n^{1-\alpha}}(1-|z|^2)^\beta|\varphi^{n-1}(z)g(z)-\psi^{n-1}(z)h(z)|\\
&=&\sup\limits_{z\in\mathbb{D}}n^\alpha(1-|z|^2)^\beta|\varphi^{n-1}(z)g(z)-\psi^{n-1}(z)h(z)|\\
\end{eqnarray*}
which implies that $\sup\limits_{n\in\mathbb{N}}\sup\limits_{z\in\mathbb{D}}(n+1)^\alpha(1-|z|^2)^\beta|\varphi^n(z)g(z)-\psi^n(z)h(z)|<\infty.$

(iv)$\Rightarrow$(iii) From proposition 3.2, it is easy to obtain.

(iii)$\Rightarrow$(ii) From proposition 3.1, it is easy to obtain.

(ii)$\Rightarrow$(i) Let $f\in \mathcal{B}^\alpha$ and $\|f\|_{\mathcal{B}^\alpha}\leq1$. By Lemma 2.1, we have
\begin{eqnarray*}
& &\|(C^g_\varphi-C^h_\psi)f\|_{\mathcal{B}^\beta}\leq\sup\limits_{z\in\mathbb{D}}|D_{\varphi,g}(z)-D_{\psi,h}(z)|(1-|\varphi(z)|^2)^\alpha|f'(\varphi(z))|\\
&+&\sup\limits_{z\in\mathbb{D}}|D_{\psi,h}(z)|\Big|(1-|\varphi(z)|^2)^\alpha f'(\varphi(z))-(1-|\psi(z)|^2)^\alpha f'(\psi(z))\Big|\\
&\lesssim&\sup\limits_{z\in\mathbb{D}}|D_{\varphi,g}(z)-D_{\psi,h}(z)|+\sup\limits_{z\in\mathbb{D}}|D_{\psi,h}(z)|\rho(z).
\end{eqnarray*}
Hence, we obtain
\begin{eqnarray*}
\|C^g_\varphi-C^h_\psi\|_{\mathcal{B}^\alpha\rightarrow\mathcal{B}^\beta}
\lesssim\sup\limits_{z\in\mathbb{D}}|D_{\varphi,g}(z)-D_{\psi,h}(z)|+\sup\limits_{z\in\mathbb{D}}|D_{\psi,h}(z)|\rho(z).
\end{eqnarray*}
Similarly
\begin{eqnarray*}
\|C^g_\varphi-C^h_\psi\|_{\mathcal{B}^\alpha\rightarrow\mathcal{B}^\beta}
\lesssim\sup\limits_{z\in\mathbb{D}}|D_{\varphi,g}(z)-D_{\psi,h}(z)|+\sup\limits_{z\in\mathbb{D}}|D_{\varphi,g}(z)|\rho(z).
\end{eqnarray*}
Thus, the statements (i),(iia),(iii),(iv) are equivalent. Similarly, the statements (i),(iib),(iii),(iv) are equivalent. The proof of the theorem is complete.
\section{Essential norm of $C^g_\varphi-C^h_\psi$: $\mathcal{B}^\alpha\rightarrow\mathcal{B}^\beta$}
In this section we give an estimate for essential norm of $C^g_\varphi-C^h_\psi$ from $\mathcal{B}^\alpha$ to $\mathcal{B}^\beta$. We need some auxiliary results. For $r\in (0,1)$, let $K_rf(z)=f(rz)$. Then $K_r$ is a compact operator on the space $\mathcal{B}^\alpha$ or $\mathcal{B}^\alpha_0$ for any positive number $\alpha$, with $\|K_r\|_{\mathcal{B}^\alpha}\leq1$. The following lemma can be found in \cite{Zh}:

{\bf Lemma 4.1.} Let $0<\alpha<\infty$. Then there is a sequence $\{r_k\}_{k=1}^\infty$ with $0<r_k<1$ tending to 1, such that the compact operator
$
L_n=\frac{1}{n}\sum\limits_{k=1}^nK_{r_k}
$
acting on $\mathcal{B}_0^\alpha$ satisfies:

(i) For any $t\in[0,1), \lim\limits_{n\rightarrow\infty}\sup\limits_{\|f\|_{\mathcal{B}^\alpha}\leq 1}\sup\limits_{|z|\leq t}|[(I-L_n)f]'(z)|=0$.

(ii) For any $s\in[0,1), \lim\limits_{n\rightarrow\infty}\sup\limits_{\|f\|_{\mathcal{B}^\alpha}\leq 1}\sup\limits_{|z|\leq s}|(I-L_n)f(z)|=0$.

(iii) $\limsup\limits_{n\rightarrow\infty}\|I-L_n\|\leq 1$.

Furthermore, these statements hold as well for the sequence of biadjoints $L_n^{**}$ on $\mathcal{B}^\alpha$.

Here we prove the following two useful propositions in this section.

{\bf Proposition 4.2.} Let $0<\alpha,\beta<\infty$. Let $\varphi, \psi\in S(\mathbb{D})$. Then the following inequalities hold:

(i) $\lim\limits_{r\rightarrow1}\sup\limits_{|\varphi(z)|>r}|D_{\varphi,g}(z)|\rho(z)\lesssim
\limsup\limits_{|a|\rightarrow1}\|(C^g_\varphi-C^h_\psi)f_a\|_{\mathcal{B}^\beta}
+\limsup\limits_{|a|\rightarrow1}\|(C^g_\varphi-C^h_\psi)g_a\|_{\mathcal{B}^\beta}.$

(ii)$\lim\limits_{r\rightarrow1}\sup\limits_{|\psi(z)|>r}|D_{\psi,h}(z)|\rho(z)\lesssim
\limsup\limits_{|a|\rightarrow1}\|(C^g_\varphi-C^h_\psi)f_a\|_{\mathcal{B}^\beta}
+\limsup\limits_{|a|\rightarrow1}\|(C^g_\varphi-C^h_\psi)g_a\|_{\mathcal{B}^\beta}.$

(iii)
\begin{eqnarray*}
\lim\limits_{r\rightarrow1}\sup\limits_{|\varphi(z)|>r\atop |\psi(z)|>r}|D_{\varphi,g}(z)-D_{\psi,h}(z)|
&\lesssim&\limsup\limits_{|a|\rightarrow1}\|(C^g_\varphi-C^h_\psi)f_a\|_{\mathcal{B}^\beta}\\
&+&\limsup\limits_{|a|\rightarrow1}\|(C^g_\varphi-C^h_\psi)g_a\|_{\mathcal{B}^\beta}.
\end{eqnarray*}

{\bf Proof.} We only need to note the proof of proposition 3.1 permit us to obtain the following consequence. For any $z\in \mathbb{D}$, we have
\begin{eqnarray*}
|D_{\varphi,g}(z)|\rho(z)
&\leq&\|(C^g_\varphi-C^h_\psi)f_a\|_{\mathcal{B}^\beta}+\|(C^g_\varphi-C^h_\psi)g_a\|_{\mathcal{B}^\beta}.
\end{eqnarray*}
Similarly
\begin{eqnarray*}
|D_{\psi,h}(z)|\rho(z)
&\leq&\|(C^g_\varphi-C^h_\psi)f_a\|_{\mathcal{B}^\beta}+\|(C^g_\varphi-C^h_\psi)g_a\|_{\mathcal{B}^\beta},
\end{eqnarray*}
and
\begin{eqnarray*}
|D_{\varphi,g}(z)-D_{\psi,h}(z)|
&\lesssim&\|(C^g_\varphi-C^h_\psi)f_{\varphi(z)}\|_{\mathcal{B}^\beta}
+\|(C^g_\varphi-C^h_\psi)f_{\psi(z)}\|_{\mathcal{B}^\beta}\\
&+&\|(C^g_\varphi-C^h_\psi)g_{\psi(z)}\|_{\mathcal{B}^\beta}.
\end{eqnarray*}

{\bf Proposition 4.3.} Let $0<\alpha,\beta<\infty$. Let $\varphi, \psi\in S(\mathbb{D})$. Then the following inequalities hold:

(i) $\limsup\limits_{|a|\rightarrow1}\|(C^g_\varphi-C^h_\psi)f_a\|_{\mathcal{B}^\beta}
\lesssim\limsup\limits_{n\rightarrow\infty}\sup\limits_{z\in\mathbb{D}}(n+1)^\alpha(1-|z|^2)^\beta|\varphi^n(z)g(z)-\psi^n(z)h(z)|.$

(ii) $\limsup\limits_{|a|\rightarrow1}\|(C^g_\varphi-C^h_\psi)g_a\|_{\mathcal{B}^\beta}
\lesssim\limsup\limits_{n\rightarrow\infty}\sup\limits_{z\in\mathbb{D}}(n+1)^\alpha(1-|z|^2)^\beta|\varphi^n(z)g(z)-\psi^n(z)h(z)|.$

{\bf Proof.} For any $a\in\mathbb{D}$ and each $N$, it follows from the proof of proposition 3.2 that
\begin {eqnarray} \label{4.1}
\|(C^g_\varphi-C^h_\psi)f_a\|_{\mathcal{B}^\beta}    \notag\
&\leq&(1-|a|^2)^\alpha|\sum_{k=0}^{N}\frac{\Gamma(k+2\alpha)}{k!\Gamma(2\alpha)}(k+1)^{-\alpha}|a|^k\\  \notag\
&\times&\sup\limits_{z\in\mathbb{D}}(k+1)^\alpha(1-|z|^2)^\beta|\varphi^k(z)g(z)-\psi^k(z)h(z)|\\
&+&(1-|a|^2)^\alpha|\sum_{k=N+1}^{\infty}\frac{\Gamma(k+2\alpha)}{k!\Gamma(2\alpha)}(k+1)^{-\alpha}|a|^k\\
&\times&\sup\limits_{n\geq N+1}\sup\limits_{z\in\mathbb{D}}(n+1)^\alpha(1-|z|^2)^\beta|\varphi^n(z)g(z)-\psi^n(z)h(z)|,   \notag\
\end {eqnarray}
Note that
\begin{eqnarray*}
\sum_{k=0}^{\infty}\frac{\Gamma(k+2\alpha)}{k!\Gamma(2\alpha)}(k+1)^{-\alpha}|a|^k\approx \frac{1}{(1-|a|^2)^\alpha},
\end{eqnarray*}
and $\sup\limits_{n\in\mathbb{N}}\sup\limits_{z\in\mathbb{D}}(n+1)^\alpha(1-|z|^2)^\beta|\varphi^n(z)g(z)-\psi^n(z)h(z)|<\infty$, let $|a|\rightarrow 1$ in (4.1) leads to
\begin{eqnarray*}
& &\limsup\limits_{|a|\rightarrow1}\|(C^g_\varphi-C^h_\psi)f_a\|_{\mathcal{B}^\beta}\\
&\lesssim&\sup\limits_{n\geq N+1}\sup\limits_{z\in\mathbb{D}}(n+1)^\alpha(1-|z|^2)^\beta|\varphi^n(z)g(z)-\psi^n(z)h(z)|.
\end{eqnarray*}
Therefore (i) holds.

Similarly for any $a\in \mathbb{D}$ and each $N$, from the proof of proposition 3.2,
\begin{eqnarray}
& &\|(C^g_\varphi-C^h_\psi)g_a\|_{\mathcal{B}^\beta}\leq \|(C^g_\varphi-C^h_\psi)f_a\|_{\mathcal{B}^\beta}\\     \notag\
&+&(1-|a|^2)^{\alpha+1}\sum_{k=1}^{\infty}
\Big(\sum_{l=0}^{k-1}\frac{\Gamma(l+2\alpha)}{l!\Gamma(2\alpha)}\Big)|a|^{k-1}\\ \notag\
&\times&\sup\limits_{z\in\mathbb{D}}(1-|z|^2)^\beta|\varphi^k(z)g(z)-\psi^k(z)h(z)|,
\end{eqnarray}
and
\begin{eqnarray*}
& &(1-|a|^2)^{\alpha+1}\sum_{k=1}^{\infty}
\Big(\sum_{l=0}^{k-1}\frac{\Gamma(l+2\alpha)}{l!\Gamma(2\alpha)}\Big)|a|^{k-1}\\
&\times&\sup\limits_{z\in\mathbb{D}}(1-|z|^2)^\beta|\varphi^k(z)g(z)-\psi^k(z)h(z)|\\
&\lesssim&(1-|a|^2)^{\alpha+1}\sum_{k=1}^{N}k^\alpha|a|^{k-1}
\sup\limits_{z\in\mathbb{D}}k^\alpha(1-|z|^2)^\beta|\varphi^k(z)g(z)-\psi^k(z)h(z)|\\
&+&(1-|a|^2)^{\alpha+1}\sum_{k=N+1}^{\infty}k^\alpha|a|^{k-1}\\
&\times&\sup\limits_{n\geq N+1}\sup\limits_{z\in\mathbb{D}}n^\alpha(1-|z|^2)^\beta|\varphi^n(z)g(z)-\psi^n(z)h(z)|\\
&\lesssim&(1-|a|^2)^{\alpha+1}\sum_{k=1}^{N}k^\alpha|a|^{k-1}
\sup\limits_{z\in\mathbb{D}}k^\alpha(1-|z|^2)^\beta|\varphi^k(z)g(z)-\psi^k(z)h(z)|\\
&+&\sup\limits_{n\geq N+1}\sup\limits_{z\in\mathbb{D}}(n+1)^\alpha(1-|z|^2)^\beta|\varphi^n(z)g(z)-\psi^n(z)h(z)|.
\end{eqnarray*}
Let $|a|\rightarrow1$ in (4.2), we get
\begin{eqnarray*}
& &\limsup\limits_{|a|\rightarrow1}\|(C^g_\varphi-C^h_\psi)g_a\|_{\mathcal{B}^\beta}\lesssim
\limsup\limits_{|a|\rightarrow1}\|(C^g_\varphi-C^h_\psi)f_a\|_{\mathcal{B}^\beta}\\
&+&\sup\limits_{n\geq N+1}\sup\limits_{z\in\mathbb{D}}(n+1)^\alpha(1-|z|^2)^\beta|\varphi^n(z)g(z)-\psi^n(z)h(z)|.
\end{eqnarray*}
Hence (ii) holds.

{\bf Theorem 4.4}  Let $0<\alpha,\beta<\infty$. Let $\varphi, \psi\in S(\mathbb{D})$. If both $C^g_\varphi$ and $C^h_\psi$ are bounded from $\mathcal{B}^\alpha$ to $\mathcal{B}^\beta$, then
\begin{eqnarray*}
\|C^g_\varphi-C^h_\psi\|_{e,\mathcal{B}^\alpha\rightarrow\mathcal{B}^\beta}
&\approx&\lim\limits_{r\rightarrow1}\sup\limits_{|\varphi(z)|>r}|D_{\varphi,g}(z)|\rho(z)
+\lim\limits_{r\rightarrow1}\sup\limits_{|\psi(z)|>r}|D_{\psi,h}(z)|\rho(z)\\
&+&\lim\limits_{r\rightarrow1}\sup\limits_{|\varphi(z)|>r\atop |\psi(z)|>r}|D_{\varphi,g}(z)-D_{\psi,h}(z)|\\
&\approx&\limsup\limits_{n\rightarrow\infty}\sup\limits_{z\in\mathbb{D}}(n+1)^\alpha(1-|z|^2)^\beta|\varphi^n(z)g(z)-\psi^n(z)h(z)|
\end{eqnarray*}

{\bf Proof.} First, we give the upper estimate. Let $\{L_n\}$ be the sequence of operators given in Lemma 4.1. Since each $L_n$ is compact as an operator from $\mathcal{B}^\alpha$ to $\mathcal{B}^\alpha$, so is $(C^g_\varphi-C^h_\psi)L_n$, and we have
\begin{eqnarray*}
\|C^g_\varphi-C^h_\psi\|_{e,\mathcal{B}^\alpha\rightarrow\mathcal{B}^\beta}
&\leq&\limsup\limits_{n\rightarrow\infty}\|C^g_\varphi-C^h_\psi-(C^g_\varphi-C^h_\psi)L_n\|_{\mathcal{B}^\alpha\rightarrow\mathcal{B}^\beta}\\
&=&\limsup\limits_{n\rightarrow\infty}\|(C^g_\varphi-C^h_\psi)(I-L_n)\|_{\mathcal{B}^\alpha\rightarrow\mathcal{B}^\beta}\\
&\leq&\limsup\limits_{n\rightarrow\infty}\sup\limits_{\|f\|_{\mathcal{B}^\alpha}\leq1}\|(C^g_\varphi-C^h_\psi)(I-L_n)f\|_{\mathcal{B}^\beta}\\
&=&\limsup\limits_{n\rightarrow\infty}\sup\limits_{\|f\|_{\mathcal{B}^\alpha}\leq1}\sup\limits_{z\in\mathbb{D}}
(1-|z|^2)^\beta\Big|[(I-L_n)f]'(\varphi(z))g(z)\\
&-&[(I-L_n)f]'(\psi(z))h(z)\Big|.
\end{eqnarray*}
For an arbitrary $r\in(0,1)$. For the sake of simplicity, we note
\begin{eqnarray*}
H^f_n(z)=(1-|z|^2)^\beta\Big|[(I-L_n)f]'(\varphi(z))g(z)-[(I-L_n)f]'(\psi(z))h(z)\Big|
\end{eqnarray*}
and set
\begin{eqnarray*}
& &\mathbb{D}_1=\{z\in\mathbb{D}:|\varphi(z)|\leq r\,\ |\psi(r)|\leq r\},\ \ \ \mathbb{D}_2=\{z\in\mathbb{D}:|\varphi(z)|\leq r\,\ |\psi(r)|>r\},\\
& &\mathbb{D}_3=\{z\in\mathbb{D}:|\varphi(z)|> r\,\ |\psi(r)|\leq r\},\ \ \ \mathbb{D}_4=\{z\in\mathbb{D}:|\varphi(z)|> r\,\ |\psi(r)|>r\}.
\end{eqnarray*}
Then
\begin{eqnarray*}
J=\sup\limits_{\|f\|_{\mathcal{B}^\alpha}\leq1}\sup\limits_{z\in\mathbb{D}}H^f_n(z)
=\max\sup\limits_{\|f\|_{\mathcal{B}^\alpha}\leq1}\sup\limits_{z\in\mathbb{D}_i}H^f_n(z)
=\max\{J_1,J_2,J_3,J_4\},
\end{eqnarray*}
where $J_i=\sup\limits_{\|f\|_{\mathcal{B}^\alpha}\leq1}\sup\limits_{z\in\mathbb{D}_i}H^f_n(z)$.
By (i) of Lemma 4.1
\begin{eqnarray*}
& &\limsup\limits_{n\rightarrow\infty}J_1
=\limsup\limits_{n\rightarrow\infty}\sup\limits_{\|f\|_{\mathcal{B}^\alpha}\leq1}\sup\limits_{z\in\mathbb{D}_1}H^f_n(z)\\
&\leq&\limsup\limits_{n\rightarrow\infty}\sup\limits_{\|f\|_{\mathcal{B}^\alpha}\leq1}\sup\limits_{|\varphi(z)|\leq r}
(1-|z|^2)^\beta|g(z)|\Big|[(I-L_n)f]'(\varphi(z))\Big|\\
&+&\limsup\limits_{n\rightarrow\infty}\sup\limits_{\|f\|_{\mathcal{B}^\alpha}}\sup\limits_{|\psi(z)|\leq r}
(1-|z|^2)^\beta|h(z)|\Big|[(I-L_n)f]'(\psi(z))\Big|\\
&=&0
\end{eqnarray*}
where we use the fact that $\sup\limits_{z\in\mathbb{D}}(1-|z|^2)^\beta|g(z)|<\infty$ and $\sup\limits_{z\in\mathbb{D}}(1-|z|^2)^\beta|h(z)|<\infty$. Since $C_\varphi^g$ and $C_\psi^h$ are bounded from $\mathcal{B}^\alpha$ to $\mathcal{B}^\beta$.
In addition,
\begin{eqnarray*}
H^f_n(z)&=&(1-|z|^2)^\beta\Big|[(I-L_n)f]'(\varphi(z))g(z)-[(I-L_n)f]'(\psi(z))h(z)\Big|\\
&\leq&|D_{\varphi,g}(z)-D_{\psi,h}(z)|(1-|\varphi(z)|^2)^\alpha\Big|[(I-L_n)f]'(\varphi(z))\Big|\\
&+&|D_{\psi,h}(z)|\Big|(1-|\varphi(z)|^2)^\alpha[(I-L_n)f]'(\varphi(z))\\
&-&(1-|\psi(z)|^2)^\alpha[(I-L_n)f]'(\psi(z))\Big|\\
&\lesssim&|D_{\varphi,g}(z)-D_{\psi,h}(z)|(1-|\varphi(z)|^2)^\alpha\Big|[(I-L_n)f]'(\varphi(z))\Big|
+|D_{\psi,h}(z)|\rho(z).
\end{eqnarray*}
Similarly
\begin{eqnarray*}
H^f_n(z)\lesssim|D_{\varphi,g}(z)-D_{\psi,h}(z)|(1-|\psi(z)|^2)^\alpha\Big|[(I-L_n)f]'(\psi(z))\Big|
+|D_{\varphi,g}(z)|\rho(z).
\end{eqnarray*}
Hence
\begin{eqnarray*}
\limsup\limits_{n\rightarrow\infty}J_2
&\lesssim&\limsup\limits_{n\rightarrow\infty}\sup\limits_{\|f\|_{\mathcal{B}^\alpha}\leq1}\sup\limits_{z\in\mathbb{D}_2}
\Big[|D_{\varphi,g}(z)-D_{\psi,h}(z)|(1-|\varphi(z)|^2)^\alpha\\
&\times&|[(I-L_n)f]'(\varphi(z))|+|D_{\psi,h}(z)|\rho(z)\Big]\\
&\lesssim&\limsup\limits_{n\rightarrow\infty}\sup\limits_{\|f\|_{\mathcal{B}^\alpha}\leq1}\sup\limits_{|\varphi(z)|\leq r}
\Big[|D_{\varphi,g}(z)-D_{\psi,h}(z)|(1-|\varphi(z)|^2)^\alpha\\
&\times&|[(I-L_n)f]'(\varphi(z))|\Big]+\sup\limits_{|\psi(z)|> r}|D_{\psi,h}(z)|\rho(z)\\
&\lesssim&\sup\limits_{|\psi(z)|> r}|D_{\psi,h}(z)|\rho(z),
\end{eqnarray*}
where we use (i) of Lemma 4.1 and the fact $\sup\limits_{z\in \mathbb{D}}|D_{\varphi,g}(z)-D_{\psi,h}(z)|<\infty$. Since $r$ is arbitrary, we have
\begin{eqnarray*}
\limsup\limits_{n\rightarrow\infty}J_2\lesssim\lim\limits_{r\rightarrow1}\sup\limits_{|\psi(z)|> r}|D_{\psi,h}(z)|\rho(z).
\end{eqnarray*}
Similarly, we can prove that
\begin{eqnarray*}
\limsup\limits_{n\rightarrow\infty}J_3\lesssim\lim\limits_{r\rightarrow1}\sup\limits_{|\varphi(z)|> r}|D_{\varphi,g}(z)|\rho(z).
\end{eqnarray*}
Also
\begin{eqnarray*}
\limsup\limits_{n\rightarrow\infty}J_4
&\lesssim&\limsup\limits_{n\rightarrow\infty}\sup\limits_{\|f\|_{\mathcal{B}^\alpha}\leq1}\sup\limits_{z\in\mathbb{D}_4}
\Big[|D_{\varphi,g}(z)-D_{\psi,h}(z)|(1-|\varphi(z)|^2)^\alpha\\
&\times&|[(I-L_n)f]'(\varphi(z))|+|D_{\psi,h}(z)|\rho(z)\Big]\\
&\lesssim&\limsup\limits_{n\rightarrow\infty}\sup\limits_{\|f\|_{\mathcal{B}^\alpha}\leq1}\sup\limits_{|\varphi(z)|>r\atop|\psi(z)|>r}
|D_{\varphi,g}(z)-D_{\psi,h}(z)|\|(I-L_n)f\|_{\mathcal{B}^\alpha}\\
&+&\sup\limits_{|\psi(z)|> r}|D_{\psi,h}(z)|\rho(z)\\
&\lesssim&\sup\limits_{|\varphi(z)|>r\atop|\psi(z)|>r}
|D_{\varphi,g}(z)-D_{\psi,h}(z)|\|(I-L_n)f\|_{\mathcal{B}^\alpha}\\
&+&\sup\limits_{|\psi(z)|> r}|D_{\psi,h}(z)|\rho(z),
\end{eqnarray*}
where we use
\begin{eqnarray*}
\limsup\limits_{n\rightarrow\infty}\|(I-L_n)f\|_{\mathcal{B}^\alpha}\leq\limsup\limits_{n\rightarrow\infty}\|(I-L_n)f\|\|f\|_{\mathcal{B}^\alpha}\leq1.
\end{eqnarray*}
Thus, we have
\begin{eqnarray*}
\limsup\limits_{n\rightarrow\infty}J_4
&\lesssim&\lim\limits_{r\rightarrow1}\sup\limits_{|\varphi(z)|>r\atop|\psi(z)|>r}
|D_{\varphi,g}(z)-D_{\psi,h}(z)|+\lim\limits_{r\rightarrow1}\sup\limits_{|\psi(z)|> r}|D_{\psi,h}(z)|\rho(z).
\end{eqnarray*}
Hence, from proposition 4.2 and 4.3, we obtain
\begin{eqnarray}
\limsup\limits_{n\rightarrow\infty}J&=&\max\Big\{\limsup\limits_{n\rightarrow\infty}J_1,\limsup\limits_{n\rightarrow\infty}J_2,  \notag\
\limsup\limits_{n\rightarrow\infty}J_3,\limsup\limits_{n\rightarrow\infty}J_4\Big\}\\           \notag\
&\lesssim&\lim\limits_{r\rightarrow1}\sup\limits_{|\varphi(z)|> r}|D_{\varphi,h}(z)|\rho(z)        \notag\
+\lim\limits_{r\rightarrow1}\sup\limits_{|\psi(z)|> r}|D_{\psi,h}(z)|\rho(z)\\                   \notag\
&+&\lim\limits_{r\rightarrow1}\sup\limits_{|\varphi(z)|>r\atop|\psi(z)|>r}|D_{\varphi,g}(z)-D_{\psi,h}(z)|\\   \notag\
&\lesssim&\limsup\limits_{|a|\rightarrow1}\|(C^g_\varphi-C^h_\psi)f_a\|_{\mathcal{B}^\beta}                      \notag\
+\limsup\limits_{|a|\rightarrow1}\|(C^g_\varphi-C^h_\psi)g_a\|_{\mathcal{B}^\beta}\\
&\lesssim&\limsup\limits_{n\rightarrow\infty}\sup\limits_{z\in\mathbb{D}}(n+1)^\alpha(1-|z|^2)^\beta|\varphi^n(z)g(z)-\psi^n(z)h(z)|
\end{eqnarray}

Next, we give the lower estimate. Let $n\geq1$. Consider the function $f_n=z^n/\|z^n\|_{\mathcal{B}^\alpha}$. Then $\|f_n\|_{\mathcal{B}^\alpha}=1$ and $f_n$ converges to 0 weakly in $\mathcal{B}^\alpha$. In particular, if $K$ is any compact operator from $\mathcal{B}^\alpha$ to $\mathcal{B}^\beta$, then $\lim\limits_{n\rightarrow\infty}\|Kf_n\|_{\mathcal{B}^\beta}=0$. Therefore
\begin{eqnarray*}
\|C_\varphi^g-C_\psi^h-K\|_{\mathcal{B}^\alpha\rightarrow\mathcal{B}^\beta}&\geq&\limsup\limits_{n\rightarrow\infty}
\|(C_\varphi^g-C_\psi^h-K)f_n\|_{\mathcal{B}^\beta}\\
&\geq&\|\limsup\limits_{n\rightarrow\infty}(C_\varphi^g-C_\psi^h)f_n\|_{\mathcal{B}^\beta}.
\end{eqnarray*}
Hence
\begin{eqnarray}
& &\|C_\varphi^g-C_\psi^h\|_{e,\mathcal{B}^\alpha\rightarrow\mathcal{B}^\beta}              \notag\
\geq\limsup\limits_{n\rightarrow\infty}(C_\varphi^g-C_\psi^h)f_n\|_{\mathcal{B}^\beta}\\   \notag\
&\gtrsim&\limsup\limits_{n\rightarrow\infty}\sup\limits_{z\in\mathbb{D}}
\frac{n}{n^{1-\alpha}}(1-|z|^2)^\beta|\varphi^{n-1}(z)g(z)-\psi^{n-1}(z)h(z)|\\ \notag\
&=&\limsup\limits_{n\rightarrow\infty}\sup\limits_{z\in\mathbb{D}}n^\alpha(1-|z|^2)^\beta|\varphi^{n-1}(z)g(z)-\psi^{n-1}(z)h(z)|\\
&=&\limsup\limits_{n\rightarrow\infty}\sup\limits_{z\in\mathbb{D}}(n+1)^\alpha(1-|z|^2)^\beta|\varphi^n(z)g(z)-\psi^n(z)h(z)|.
\end{eqnarray}
Combining (4.3) and (4.4), we immediately get the desired result. The proof is complete.

From Theorem 4.4, we immediately get the following corollary.

{\bf Corollary 4.5}  Let $0<\alpha,\beta<\infty$. Let $\varphi, \psi\in S(\mathbb{D})$. If both $C^g_\varphi$ and $C^h_\psi$ are bounded from $\mathcal{B}^\alpha$ to $\mathcal{B}^\beta$, then the following statements are equivalent.

(i) $C^g_\varphi-C^h_\psi:\mathcal{B}^\alpha\rightarrow\mathcal{B}^\beta$ is compact.

(ii)
\begin{eqnarray*}
\lim\limits_{r\rightarrow1}\sup\limits_{|\varphi(z)|> r}|D_{\varphi,h}(z)|\rho(z)
&=&\lim\limits_{r\rightarrow1}\sup\limits_{|\psi(z)|> r}|D_{\psi,h}(z)|\rho(z)\\
&=&\lim\limits_{r\rightarrow1}\sup\limits_{|\varphi(z)|>r\atop|\psi(z)|>r}|D_{\varphi,g}(z)-D_{\psi,h}(z)|=0\\
\end{eqnarray*}

(iii)
$\limsup\limits_{|a|\rightarrow1}\|(C^g_\varphi-C^h_\psi)f_a\|_{\mathcal{B}^\beta}
=\limsup\limits_{|a|\rightarrow1}\|(C^g_\varphi-C^h_\psi)g_a\|_{\mathcal{B}^\beta}=0$.

(iv)
$\limsup\limits_{n\rightarrow\infty}\sup\limits_{z\in\mathbb{D}}(n+1)^\alpha(1-|z|^2)^\beta|\varphi^n(z)g(z)-\psi^n(z)h(z)|=0$.

%\end{enumerate}

\end{document}